\documentclass[11pt,a4paper,notitlepage,final,onecolumn,oneside]{article}
\usepackage{exscale}
\usepackage{latexsym}
\usepackage{calc}
\usepackage{amssymb}
\usepackage{url}

\newcommand{\proof}{{\bfseries Proof:\quad}}
\newcommand{\p}{\ensuremath{\hat{\mathrm{p}}}}

\begin{document}

\newtheorem{lemma}{Lemma}
\newtheorem{proposition}[lemma]{Proposition}
\newtheorem{theorem}[lemma]{Theorem}
\newtheorem{corollary}[lemma]{Corollary}


\title{The $k$-orbit theory and polycirculant conjecture}
\author{Aleksandr Golubchik\\
        \small Osnabrueck, Germany,\\}

\date{15.03.2005}
\maketitle

\begin{abstract}
The paper contains a proof of the conjecture of M. Klin and D.
Maru$\breve{\rm s}$i$\breve{\rm c}$ that an automorphism group of a
transitive graph contains a permutation, decomposed in cycles of the same
length. The proof is based on the $k$-orbit theory developed by author.
\bigskip

\emph{Key words:} $k$-orbits, partitions, permutations, symmetry, groups

\end{abstract}


\section{Introduction}

A permutation group $G(V)$ is called $2$-closed if it is an automorphism
group of a graph, another if it is an automorphism group of an ordered
partition $\overrightarrow Q_2$ of a set $V^2$. If $G'<G$ and
$G'\overrightarrow Q_2=G'\overrightarrow Q_2=\overrightarrow Q_2$, then
one says that $G$ is a $2$-closure for $G'$. A permutation is called
regular if it is decomposed in cycles of the same length.

In \cite{Cameron} P. Cameron described the conjecture of M. Klin and D.
Maru$\breve{\rm s}$i$\breve{\rm c}$ (Polycirculant conjecture) that was
announced on the British Combinatorial Conference (Problem BCC15.12
(DM282)) in July 1995:

\begin{theorem}\label{PC}
Every $2$-closed transitive group contains a regular element.

\end{theorem}

Below we give a proof of this statement, based on the $k$-orbit theory
developed by author. The first attempt to describe this theory was
launched in \cite{Golubchik}.


\section{$k$-Orbit theory}

Let $G(V)$ be a permutation group of a degree $n$ and $V^{(k)}$ ($k\leq
n$) be a non-diagonal part of $V^k$ (i.e. all $k$ values of coordinates of
a $k$-tuple $\alpha_k\in V^{(k)}$ are different), then the action of $G$
on $V^{(k)}$ forms a partition $GV^{(k)}$ of $V^{(k)}$ on $G$-invariant
classes. This partition is called a system of $k$-orbits of $G$ and we
denote it as $Orb_k(G)$.

The $k$-orbit theory studies symmetry properties of $k$-orbits that cannot
be obtained from permutation group theory, studying a permutation group
$G$ as a permutation algebra and therefore being not able to see the
inside structure of a $n$-orbit of $G$ and relations between $k$-orbits of
subgroups of $G$. Also the graph representation of permutation groups that
considers an action of a permutation group on $V^2$ is not able to see
all relations in $Orb_k(G)$ for $2<k\leq n$.

The $k$-orbit theory gives a new view on abstract finite groups and finite
permutation groups and shows a way to simple solutions of some problems
which either are not solved or have a complicated solution. These problems
are: the polynomial solution of graphs isomorphism problem; the problem of
a full invariant of finite groups, admitting a primitive (isomorphic)
permutation representation; a simple proof of the Feit-Thompson theorem:
the solvability of finite groups of odd order (that lies in the foundation
of the Classification of finite simple groups \cite{Gorenstein}); The
simple proof of the Fein-Kantor-Schacher theorem: any transitive
permutation group contains a fixed-point-free prime-power element
\cite{FKS}; and the problem that we consider in this paper.

The specificity of the $n$-orbit representation of a group $G$ is a
possibility to do a group visible. In order to go to this visibility one
makes a partition of a matrix of a $n$-orbit on cells of $k$-orbits of
subgroups of an investigated permutation group $G$ and studies symmetry
properties of cells and the whole partition. Below we consider some facts
from the $k$-orbit theory which are necessary for the proof of theorem
\ref{PC}.

The $k$-orbit theory considers different actions of permutations on
$k$-sets and symmetries that follows from those actions. Let
$\alpha_n=\langle v_1\ldots v_n\rangle$ be a $n$-tuple, then the left
action of a permutation $g$ on $\alpha_n$ is defined as $g\alpha_n=\langle
gv_1\ldots gv_n\rangle$ and the right action of $g$ on $\alpha_n$ is given
as $\alpha_ng=\langle v_{g1}\ldots v_{gn}\rangle$ (here $V=[1,n]$). From
that we obtain the left and right action of a permutation on $k$-tuples
and $k$-sets.

The left action of a permutation is (on definition) an isomorphism
(\,$(123)\{123,132\}=\{312,321\}$\,), the right action is not an
isomorphism (\,$\{123,132\}(123)=\{312,213\}$\,). If the right action is an
isomorphism, then it is an indicium for existing of a subgroup with
non-trivial normalizer, whose $n$-orbit contains (right-)isomorphic
$k$-subsets connected with considered right action.

Let $X_n$ be a $n$-orbit of $G$. One of $n$-tuples from $X_n$ we choose as
an \emph{initial} $n$-tuple and consider all permutations from $G$
relatively to this $n$-tuple. If $\alpha_n$ is an initial $n$-tuple, then
$X_n=G\alpha_n$. The specificity of an initial $n$-tuple is the equality
of number and order values of its coordinates (that we take from the same
ordered set $V$). We call the initial $n$-tuple as $n$-space and its
ordered subsets as subspaces. A $k$-orbit $X_k$ of $G$ is a projection of
$X_n$ on some subspace $I_k\in X_k$. We write this projection as
$X_k=\p^{(I_k)}X_n$. So we consider subspaces as ordered sets of
coordinates and $k$-orbits as non-ordered sets of $k$-tuples. The
different orderings of $k$-tuples (lines) in a matrix of a $k$-orbit shows
different symmetry properties of $k$-orbit, related with corresponding
properties of the investigated permutation group.

We say that two $k$-sets $X_k$ and $X_k'$ are $G$-isomorphic if they are
isomorphic and connected with a permutation $g\in G$:  $X_k'=gX_k$. If we
study a group $G$ and $X_k,X_k'$ are isomorphic, but not $G$-isomorphic,
then we say that they are $S_n$-isomorphic (where $S_n$ is the symmetric
group).

We define the left action $GY_n$ of a group $G$ on a $n$-subset
$Y_n$ as $GY_n\equiv \{gY_n:\, g\in G\}$. The left action $AX_n$ of a
subgroup $A<G$ on a $n$-orbit $X_n$ of $G$ is given as $AX_n\equiv
\{A\alpha_n:  \alpha_n\in X_n\}$. Correspondingly the right action
$Y_nG\equiv \{Y_ng:\, g\in G\}$ and $X_nA\equiv \{\alpha_nA:\alpha_n\in
X_n\}$.  The left action of $G$ on a $k$-set $Y_k$ and $A$ on a $k$-orbit
$X_k$ of $G$ for $k<n$ is the same as for $k=n$, because
$\p^{(I_k)}g\alpha_n=g\p^{(I_k)}\alpha_n$.  By the right action the
equality $\p^{(I_k)}(\alpha_ng)= (\p^{(I_k)}\alpha_n)g$ has sense only if
$\alpha_n$ is the initial $n$-tuple. Hence the right action on $k$-sets
has not a direct reduction from the right  action on $n$ sets. If $X_n$ is
a $n$-orbit of $G$, $A<G$ and $Y_n$ is a $n$-orbit of $A$, then
$L_n=GY_n=X_nA$ and $R_n=Y_nG=AX_n$ are partitions of $X_n$ on $n$-orbits
of left and right cosets of $A$ in $G$, at that $n$-orbits of left cosets
are $n$-orbits of conjugate with $A$ subgroups of $G$ and $n$-orbits of
right cosets are $n$-orbits of $A$.  Projections of $X_n$, $L_n$ and $R_n$
on a subspace $I_k$ gives a covering $L_k$ and a partition $R_k$ of the
$k$-orbit $X_k$ on $k$-orbits of left and right cosets of $A$ in $G$.  If
$Y_k=\p^{(I_k)}Y_n$, then $L_k=GY_k$ and $R_k=AX_k$. Below under action of
permutations we understand the left action.

Let $k$ be a divisor of $|G|$. We call $k$ as an \emph{automorphic number}
if there exists a $k$-element \emph{suborbit} of $G$ (a $k$-element orbit
of a subgroup of $G$).

Let $Co(\alpha_k)$ be a set of coordinates of a $k$-tuple $\alpha_k$, then
for a $k$-set $X_k$ we define $Co(X_k)$ as $\{Co(\alpha_k):\,\alpha_k\in
X_k\}$.

If $Co(\alpha_k)$ is a suborbit of $G$, then we say that $\alpha_k$ is an
\emph{automorphic $k$-tuple} and $Co(\alpha_k)$ is \emph{an automorphic
subset} of $V$.  Let $X_k$ be a $k$-set we say that $X_k$ is an
\emph{automorphic $k$-set}, if $Aut(X_k)$ acts on $X_k$ transitive. Let
$X_k$ be a $k$-orbit (an automorphic $k$-set) and $\alpha_k\in X_k$ be an
automorphic $k$-tuple, then a $k$-orbit $X_k$ we call as
\emph{right-automorphic $k$-orbit} or \emph{$k$-rorbit}.

Let $U$ be a set and $Q$ be a set of subsets of $U$, then $\cup Q\equiv
\cup_{(W\in Q)}W$. Under $\sqcup Q$ we understand a set of unions of
intersected on $U$ classes of $Q$. So $\sqcup Q$ is a partition of $\cup
Q$ (possibly trivial). A union and intersection of partitions $P$ and $R$
of $U$ we write as $P\sqcup R$ and $P\sqcap R$. Let $P'$ be a subpartition
of $P$, then we write $P'\sqsubset P$.

Let $\alpha_k\in X_k$ and $Y_k$ be a maximal subset of $X_k$ with
$Co(Y_k)=Co(\alpha_k)$, then we call $Y_k$ as a \emph{$k$-block} of $X_k$.

Let $\cup Co(X_k)=V$ and $\sqcup Co(X_k)$ be a non-trivial partition of
$V$, then we say that $X_k$ is \emph{incoherent}. If $\sqcup Co(X_k)$ is a
non-trivial covering of $V$, then we say that $X_k$ is \emph{coherent}. Any
coherent $k$-orbit consists of intersected on $V$ $k$-blocks. Classes of
a partition of an incoherent $k$-orbit are coherent $k$-suborbits or
$k$-blocks. A coherent $k$-orbit $X_k$ is defined on a set $V=\cup
Co(X_k)$. Its coherent or incoherent $k$-suborbit $Y_k$ is defined on a
set $U=\cup Co(Y_k)\subset V$. In order to show this we say that $Y_k$ is
$U$-coherent ($U$-incoherent). If a ($V$-) coherent $k$-orbit contains no
$U$-coherent and no $U$-incoherent $k$-suborbit for any $U\subset V$, then
we call it as \emph{elementary coherent}.

Let $Y_k$ be a $k$-suborbit of $G$. The maximal transitive on $Y_k$
subgroup of $G$ we call a \emph{stabilizer} of $Y_k$ in $G$ and write it as
$Stab_{\, G}(Y_k)$ or simply $Stab(Y_k)$.

Let $X_k$ be a $k$-orbit of $G$ and $Y_k\subset X_k$ be a $k$-block, then
evidently $L_k=GY_k$ is a partition of $X_k$ and $Y_k$ is a $k$-orbit of
$Stab(Y_k)=Stab(Co(Y_k))$.

Let a group $G$ be imprimitive, then $G$ contains a non-trivial ($1<k<n$)
incoherent $k$-orbit. If $G$ is primitive, then any (non-trivial)
$k$-orbit of $G$ is coherent. From here follows that it is convenient to
consider primitive Abelian groups as trivial imprimitive or trivial
primitive, because, as distinct from non-Abelian primitive groups, such a
group $G$ contains no non-trivial suborbit $U\subset V$ that forms a
covering $GU$ of $V$. So further under a primitive group we understand a
non-Abelian primitive group. Let a group $G$ be imprimitive and $Q$ be a
$G$-invariant partition of $V$ ($GQ=Q$), then we say that classes of $Q$
are \emph{imprimitivity blocks}.

We say that a transitive group $G(V)$ is a \emph{minimal degree group} or
\emph{md-group} if it cannot be represented isomorphically on a partition
$Q$ of $V$, otherwise we say that $G$ is a \emph{non-minimal degree group}
or \emph{nmd-group}. It follows that primitive groups are md-groups. Let
$F$ be a finite group, $A<F$ and a representation $F(F/A)$ be a md-group
isomorphic to $F$, then we say that $A$ is a md-stabilizer of $F$. So $A$
is a maximal by including subgroup of $F$ containing no normal subgroup of
$F$.

Let $\alpha_k=\langle v_1\ldots v_k\rangle$ and $\alpha_l'=\langle
v_1'\ldots v_l'\rangle$, then a $(k+l)$-tuple $\alpha_{(k+l)}=\langle
v_1\ldots v_kv_1'\ldots v_l'\rangle$ we call a \emph{concatenation} of
$\alpha_k$ and $\alpha_l'$ and write this as $\alpha_{(k+l)}=\alpha_k\circ
\alpha_l'$.

\begin{proposition}\label{Stab(Co(a_k))}
Let $\alpha_k$ be an automorphic $k$-tuple, $G=Stab(Co(\alpha_k))$,
$H=Stab(\alpha_k)$ and $X_k=G\alpha_k$, then $H$ is a normal subgroup
of $G$ and a factor group $G/H$ is isomorphic to $Aut(X_k)$.

\end{proposition}
\proof
If $H$ is not trivial, then $G$ is an intransitive group. Let $X_n$ be a
$n$-orbit of $G$, then there exists only one partition $P_n$ of $X_n$ on
$|G|/|H|$ classes with projection $\p^{(\alpha_k)}P_n=X_k$. Hence
$X_nH=HX_n=P_n$. The group $Aut(X_k)$ is isomorphic to the factor group
$G/H$, because $|Aut(X_k)|=|G/H|$ and it acts on cosets of $H$ in $G$
transitive. $\Box$\bigskip

\begin{proposition}\label{L_k=R_k}
Let $X_k\in Orb_k(G)$, $Y_k\subset X_k$ be a $k$-suborbit of $G$,
$H=Stab(Y_k)$, $L_k=GY_k$ and $R_k=HX_k$. Let $L_k=R_k$, then $H$ is a
normal subgroup of $G$.

\end{proposition}
\proof
There exists only one partition $P_n$ of $X_n$ on $|G|/|H|$ classes so
that $\p^{(I_k)}P_n=L_k$. Hence $X_nH=HX_n=P_n$. $\Box$\bigskip

\begin{proposition}\label{impr.norm}
Let $G(V)$ be an imprimitive group, $Q$ be a partition of $V$ on
imprimitivity blocks and $\overrightarrow{Q}$ be ordered set $Q$, then
$H=Stab(\overrightarrow{Q})$ is a normal subgroup of $G$.

\end{proposition}
\proof
Subgroups $\{Stab(U)<G:\, U\in Q\}$ form a class of conjugate subgroups
and $H=\cap_{(U\in Q)}Stab(U)$. If $G(V)$ is a nmd-group and $G(Q)$ is an
isomorphism, then $H$ is trivial. $\Box$\bigskip

\begin{proposition}\label{cap,cup}
Let $X_k\in Orb_k(G)$, $Y_k,Z_k\subset X_k$ be $k$-suborbits of $G$,
$L_k'=GY_k$ and $L_k''=GZ_k$, then $L_k^1=L_k'\sqcap L_k''$ and
$L_k^2=L_k'\sqcup L_k''$ are also partitions of $X_k$ on isomorphic
$k$-suborbits of $G$. If $Y_k\cap Z_k=T_k$ be not trivial and $U_k\in
L_k^2$ contains $Y_k$ and $Z_k$, then $L_k^1=GT_k$, $L_k^2=GU_k$,
$Stab(T_k)=Stab(Y_k)\cap Stab(Z_k)$ and
$Stab(U_k)=gr(Stab(Y_k),Stab(Z_k))$.

\end{proposition}
\proof
Let $k=n$, then equalities follows from corresponding properties of sets of
left cosets of subgroups of $G$. For $k<n$ equalities are projections of
corresponding equalities for $k=n$. $\Box$\bigskip

\begin{lemma}\label{|Aut(X_k)|=|X_k|}
Let for every $k$-suborbit $Y_k$ of a $k$-orbit $X_k$ a set
$L_k=Aut(X_k)Y_k$ be a partition of $X_k$, then $|Aut(X_k)|=|X_k|$.

\end{lemma}
\proof
Let $G=Aut(X_k)$, $M$ be a set of partitions of $X_k$ on isomorphic
$k$-suborbits of $G$ and $L_k',L_k''\in M$, then $L_k^1=L_k'\sqcap L_k''$
and $L_k^2=L_k'\sqcup L_k''$ are partitions from $M$. Let $L_k''\in M$ be
a partition of $X_k$ on $k$-blocks, then there exists a partition
$L_k'\in M$ so that $L_k^2=L_k'\sqcup L_k''$ is not trivial, i.e.
$|L_k^2|<min(|L_k|,|B_k|)$. By repeating we obtain that there exist
partitions $L_k',L_k''\in M$ so that $L_k^2=L_k'\sqcup L_k''=\{X_k\}$.

For each partition $L_k$ in the union process there exists a subgroup
$A<G$ so that a class $Y_k\in L_k$ is a $k$-orbit of $A$ and $|A|=|Y_k|$,
because we can begin from a subgroup that is isomorphic to the
automorphism group of a $k$-block, and with a subgroup of a prime order
that permutes $k$-blocks. Thus $G$ contains a subgroup $H$ of order
$|X_k|$ that acts transitive on $X_k$.

Let $g$ be any permutation from $G$ and $R_k=gr(g)X_k$, then $(gr(g)\cap
B)X_k=R_k$, because any class $Z_k$ of $R_k$ generates a partition
$GZ_k\in M$. It follows that the action of $G$ and $B$ on $X_k$ are
isomorphic and hence $G=B$. $\Box$\bigskip

\begin{proposition}\label{X_k.incoh}
Let $X_k$ be an incoherent automorphic $k$-set, then $|Aut(X_k)|/|X_k|>1$.

\end{proposition}
\proof
If $X_k$ is incoherent, then $Aut(X_k)$ has evidently a permutation that
fixes a $k$-tuple from $X_k$. $\Box$\bigskip

\begin{lemma}\label{Aut(el.coh)}
Let $X_k$ be an elementary coherent $k$-orbit, then $|Aut(X_k)|=|X_k|$.

\end{lemma}
\proof
Let $Y_k\subset X_k$ be a $k$-block and $Z_k\subset X_k$ be a $k$-suborbit
of $Aut(X_k)$ that has non-trivial intersection with $Y_k$. Let
$A=gr(Stab(Y_k),Stab(Z_k))$, then $k$-orbit of $A$ in $X_k$ is $X_k$,
because the elementary coherent $k$-orbit $X_k$ is the unique
super-$k$-suborbit for any $k$-block. Hence $Aut(X_k)Z_k$ is a partition
of $X_k$. Then the statement follows from lemma \ref{|Aut(X_k)|=|X_k|}.
$\Box$\bigskip

\begin{lemma}\label{a_k-circ-a_k'}
Let $G$ be a transitive group, a $n$-tuple $\alpha_n=\alpha_k\circ\,
\alpha_k'\circ\,\alpha_m$, $2k+m=n$, $\alpha_k$ and $\alpha_k'$ be
automorphic and $G$-isomorphic, and $Co(\alpha_m)$ contains no $k$-tuple
$G$-isomorphic to $\alpha_k$, then
\begin{itemize}
\item
the $(2k)$-tuple $\alpha_{2k}=\alpha_k\circ\,\alpha_k'$ is automorphic,
\item
the group $A=Stab(Co(\alpha_k))\cap Stab(Co(\alpha_k'))$ is not trivial,

\item
the group $A$ has a normalizer $N$ that acts imprimitive on
$Co(\alpha_{2k})$ with imprimitivity blocks $Co(\alpha_k)$ and
$Co(\alpha_k')$.

\end{itemize}
\end{lemma}
\proof
Let $X_n\in Orb_n(G)$, $X_k=\p^{(\alpha_k)}X_n$,
$X_k'=\p^{(\alpha_k')}X_n$ and $X_{2k}=\p^{(\alpha_{2k})}X_n$.
Let $Y_k,Y_k'\subset X_k=X_k'$ be $k$-blocks with coordinates
$Co(Y_k)=Co(\alpha_k)$ and $Co(Y_k')=Co(\alpha_k')$. Consider a maximal
$(2k)$-subset $Y_{2k}\subset X_{2k}$, whose projection on $\alpha_k$ is
$Y_k$ and a maximal $(2k)$-subset $Y_{2k}'\subset X_{2k}$, whose
projection on $\alpha_k'$ is $Y_k'$. The $(2k)$-subset $Y_{2k}$ can be
partition on $(2k)$-orbits of $Stab(Co(\alpha_k))$ and the $(2k)$-subset
$Y_{2k}'$ can be partition on $(2k)$-orbits of $Stab(Co(\alpha_k'))$.
Take in opinion that $\alpha_k$ and $\alpha_k'$ have on condition
non-trivial intersection with any $k$-tuple from $X_k\setminus
Y_k\setminus Y_k'$ and that the intersection $Y_{2k}\cap Y_{2k}'$ contains
$\alpha_{2k}$, we obtain that $Z_{2k}=Y_{2k}\cap Y_{2k}'$ is a
$(2k)$-orbit of the non-trivial subgroup $A$.

Let now $T_{2k},T_{2k}'\subset X_{2k}$ be maximal $(2k)$-subsets with
projection $Y_k\cup Y_k'$ on $\alpha_k$ and $\alpha_k'$ correspondingly,
then $U_{2k}=T_{2k}\cap T_{2k}'$ is a $(2k)$-orbit of the imprimitive on
$\alpha_{2k}$ normalizer $N$ of $A$ with imprimitivity blocks
$Co(\alpha_k)$ and $Co(\alpha_k')$. $\Box$\bigskip

For the greater than two automorphic $G$-isomorphic subsets we have more
possibilities for permutations of subsets and therefore find more possible
properties.

\begin{lemma}\label{a_k-circ-...}
Let $\alpha_n=\alpha_k^1\circ\,\ldots\circ\,\alpha_k^l\circ\,\alpha_m$,
$r=lk$, $n=r+m$, all $\alpha_k^i$ are automorphic and $G$-isomorphic and
no class of $GCo(\alpha_k^1)$ belongs to $Co(\alpha_m)$, then $r$-tuple
$\alpha_{r}=\alpha_k^1\circ\,\ldots\circ\,\alpha_k^l$ is automorphic, a
subgroup $A=\cap_{(i=1,l)} Stab(Co(\alpha_k^i))$ can be trivial and has a
normalizer $N$ that acts imprimitive on $\alpha_{r}$ with imprimitivity
blocks $Co(\alpha_k^i)$.

\end{lemma}
\proof
Let $X_k=\p^{(\alpha_k^i)}X_n$ and $X_r=\p^{(\alpha_r)}X_n$. If
$g\in G$ maps $\alpha_k^1$ on $\alpha_k^2$, then $\alpha_k^2$ can be mapped
on any $k$-tuple with coordinates $Co(\alpha_k^i)$. So in this case $A$
can be trivial, but $N$ is not trivial, because subsets $Co(\alpha_k^i)$
are non-intersected in pairs and, on condition, they are intersected with
any $k$-tuple from $X_k\setminus \cup_{(i=1,l)}Y_k^i$, where $Y_k^i$ are
$k$-blocks with $Co(Y_k^i)=Co(\alpha_k^i)$. $\Box$\bigskip

\begin{theorem}\label{Q}
Let $G(V)$ be a primitive group and $k$ be a maximal automorphic divisor of
$n$, then there exists a partition $Q$ of $V$ on automorphic
$G$-isomorphic $k$-subsets or on automorphic $S_n$-isomorphic $k$-subsets.

\end{theorem}
\proof
Let the statement be not correct, then there exists a subset $U$ of $V$
that can be partition on automorphic $G$-isomorphic $k$-tuples
$\{\alpha_k,\ldots\}=P$ (\,where $1<|P|<n/k$\,), and the subset
$W=V\setminus U$ contains no automorphic $k$-tuple $G$-isomorphic to
$\alpha_k$, then from lemma \ref{a_k-circ-...} immediately follows that
$U$ is an automorphic subset.  If $|U|$ divides $n$, then $k$ is not the
maximal divisor of $n$.  Contradiction. So $|U|$ does not divide $n$ and
$Q'=GU$ is a covering of $V$. Let $Q''=GCo(\alpha_k)$, then $Q'$ and $Q''$
are $G$-invariant and therefore an intersection of any two classes of $Q'$
contains whole number of $k$-tuples from $Q''$. Hence $W$ must contain
some class from $Q''$.  Contradiction.~$\Box$\bigskip

\begin{lemma}\label{nmd-2.cl}
A nmd-group is $2$-closed.

\end{lemma}
\proof
Let $G(V)$ be not $2$-closed a nmd-group and $A(V)$ be a $2$-closure of
$G(V)$, then $A(V)$ and $G(V)$ have the same set of imprimitivity
blocks $Q$ and hence $A(V)$ is a md-group. Let $G(Q)$ be a md-reduction
of $G(V)$ and $A(Q)$ be a md-reduction of $A(V)$. Since $G(V)$ is not
$2$-closed, then $G(Q)$ is not $2$-closed too and hence $A(Q)\neq G(Q)$.
It follows that $A(Q)$ is isomorphic to $A(V)$, i.e. $A(V)$ is a
nmd-group. Contradiction. $\Box$\bigskip

\begin{lemma}\label{nmd->Xk=Xn}
Let $Q$ be a set of imprimitivity blocks of a nmd-group $G(V)$ and
$Co(I_k)\in Q$, then $Stab(I_k)$ is trivial.

\end{lemma}
\proof
Let $Q=\{Co(I_k^i):\,i\in [1,|Q|]\}$, $X_n\in Orb_n(G)$,
$X_k=\p^{(I_k^i)}X_n$ and $X_{2k}^{ij}=\p^{(I_k^i\circ I_k^j)}X_n$, then
from lemma \ref{nmd-2.cl} it follows that $G=\cap_{(i,j\in
[1,|Q|])}Aut(X_{2k}^{ij})$. Let $|X_n|/|X_k|>1$, then there exists a
partition $L_{2k}^{ij}$ of $X_{2k}^{ij}$ so that
$\p^{(I_k^i)}L_{2k}^{ij}=\p^{(I_k^j)}L_{2k}^{ij}= \{X_k\}$. It follows
that $Aut(X_{2k}^{ij})$ and hence $G$ contain a subgroup $H$ of order
$|X_k|$. Then there exists only one partition $P_n$ of $X_n$ on
$|X_n|/|X_k|$ automorphic classes so that $\p^{(I_k^i)}P_n= \{X_k\}$,
hence $H$ is a normal subgroup of $G$. Let a representation $G(Q)$ be a
md-group, then $Stab(Co(I_k))$ is a md-stabilizer of $G$. But
$Stab(Co(I_k))$ has non-trivial intersection with $H$ and hence it cannot
be a md-stabilizer of $G$. Contradiction. $\Box$\bigskip

\begin{lemma}\label{nmd-prime}
Let $G(V)$ be a nmd-group, then there exists an isomorphic representation
$G(Q)$ of $G(V)$ on imprimitivity blocks of a prime power.

\end{lemma}
\proof
Let $|V|/|Q|=k$ be not prime and $p$ be a prime divisor of $k$. Let
$Co(I_k)\in Q$, then $I_k$ contains an automorphic $p$-tuple $I_p$, so
that $Q'=GCo(I_p)\sqsubset GCo(I_k)=Q$. Hence $G$ has an imprimitive
representation on classes of $Q'$. If $G(Q')$ is a homomorphism, then
$G(Q)$ is a homomorphism too.  Contradiction. Thus $G(Q')$ is an
isomorphism.~$\Box$\bigskip

This statement can be generalized:

\begin{lemma}\label{md-prime}
Let $Q$ be a partition of $V$ and $G(Q)$ be a homomorphism of $G(V)$,
then there exists a partition $Q'\sqsubset Q$ on subsets of $V$ of a prime
power so that $G(Q')$ is also a homomorphism, at that if $G(Q)$ is an
isomorphism, then $G(Q')$ is an isomorphism too.

\end{lemma}


\section{The proof of theorem \ref{PC}}

\subsection{Imprimitive nmd-groups}

From lemma \ref{nmd-prime} it follows that any nmd-group contains a
regular element.

\subsection{Primitive groups}

Let $G$ be a primitive group. If $G$ contains a nmd-subgroup, then it
contains a regular element. Let $G$ contains no nmd-subgroup, but an
imprimitive md-subgroup $A$.

\begin{lemma}\label{md-impr<pr}
Let $A$ be a maximal imprimitive md-subgroup of $G$, then $G$ is
$2$-closed if and only if $A$ is $2$-closed.

\end{lemma}
\proof
Let $Q$ be a set of imprimitivity blocks of $A$, then $A=Stab(Q)$ and all
classes of $GQ$ (partitions of $V$) are different. Let $X_n\in Orb_n(G)$,
$Y_n\in Orb_n(A)$ and $L_n=GY_n$. Let $A$ be not $2$-closed, then any
class $Z_n\in L_n$ is not $2$-closed and hence $X_n$ is not $2$-closed.
Let $X_n$ be not $2$-closed and $X_n'$ be a $2$-closure of $X_n$, then
$X_n'$ contains an imprimitive on $Q$ $n$-suborbit $Y_n'$ that is a
$2$-closure of $Y_n$. $\Box$\bigskip\\
Thus a consideration of a primitive group $G$ can be reduced to a
consideration of its imprimitive md-subgroup. So we have to consider a
primitive group that contains no imprimitive subgroup.

\begin{lemma}\label{no.impr.subgr}
Let $G$ contains no imprimitive subgroup, then it is not $2$-closed.

\end{lemma}
\proof
If $G$ contains no imprimitive subgroup, then (lemma \ref{a_k-circ-...}
and theorem \ref{Q}) there exists a partition $Q$ of $V$ on automorphic
$S_n$-isomorphic $k$-subsets and no partition of $V$ on automorphic
$G$-isomorphic $k$-subsets. Let $X_n\in Orb_n(G)$,
$Q=\{Co(\alpha_k^i):\,i\in [1,|Q|])\}$, $X_k^i=\p^{(\alpha_k^i)}X_n$ and
$L_k=\{X_k^i:\,i\in [1,|Q|])\}$, then any two $k$-blocks from a $k$-orbit
$X_k^i\in L_k$ are intersected on $V$. Hence a non-trivial normalizer
$N=Stab_{\,S_n}L_k$ of $G$ is a $2$-closure of $G$. $\Box$\bigskip

So we can reduce the consideration of theorem \ref{PC} to imprimitive
md-groups.

\subsection{Imprimitive md-groups}

Let $G(V)$ be an imprimitive md-group and $Q$ be a partition of $V$ on
imprimitivity blocks, then a group $H=Stab(\overrightarrow{Q})$ is not
trivial. Let $Y_n\in Orb_n(H)$, $Q=\{Co(I_k^i):\,i\in [1,|Q|]\}$,
$Y_{2k}^{ij}=\p^{(I_k^i\circ I_k^j)}Y_n$ and $G$ be $2$-closed, then $H$
is also $2$-closed and hence $H=\cap_{(i,j\in
[1,|Q|])}Aut(Y_{2k}^{ij};V)$, where $Aut(Y_{2k}^{ij};V)$ is an extension
of $Aut(Y_{2k}^{ij})$ on $V$. Accordingly to lemma \ref{md-prime} we can
assume that $|V|/|Q|=p$ is prime, then $Aut(Y_{2k}^{ij};V)$ (for any
$i,j\in [1,|Q|]$) is a direct product and contains a regular element of
order $p$. Hence $H$ contains a regular element of order $p$.


\section*{Conclusion}

The main inference from the investigation of $k$-orbits symmetry
properties is that a finite group $F$ is not closed by its algebraic
properties, because the group algebra $F$ is equivalent to the action of
$F$ on $F$, but this algebra generates also the action of $F$ on $F^k$.
Properties of that action not always can be interpreted with group algebra
or with traditional permutation group theory. Namely such properties are
the subject of investigation of the $k$-orbit theory. Some investigations
related with an application of $k$-orbit representation to problems
announced above are described in \cite{graph.iso} - \cite{FinGr}

The $k$-orbit theory gives a new view on finite many-dimensional
symmetries and can bring new ideas and applications to other sciences
studying symmetrical objects.


\end{document}